\documentclass[11pt,reqno]{amsart}

\usepackage{euscript}
\usepackage{epsfig}
\usepackage{amssymb}

\theoremstyle{plain}

\theoremstyle{definition}

\newcommand{\bp}{{\mathbb P}}
\newcommand{\bz}{{\mathbb Z}}
\newcommand{\bq}{{\mathbb Q}}
\newcommand{\bn}{{\mathbb N}}
\newcommand{\cc}{\bar{C}}
\newcommand{\vp}{\varphi}

\begin{document}

\title{Classification of rational unicuspidal projective curves
whose singularities have one Puiseux pair}
\author{J. Fern\'andez de Bobadilla,
I. Luengo, A. Melle-Hern\'andez, and A. N\'emethi}
\address{Department of Mathematics, University of Utrecht, Postbus 80010, 3508TA Utrecht, The Netherlands}
\address{Facultad de Matem\'aticas\\ Universidad Complutense\\ Plaza de Ciencias\\ E-28040, Madrid, Spain}
\address{Department of Mathematics\\Ohio State University\\Columbus,
OH 43210,USA; and R\'enyi Institute of Mathematics, Budapest,
Hungary.}
\address{}
\email{bobadilla@math.uu.nl} \email{iluengo@mat.ucm.es}
\email{amelle@mat.ucm.es} \email{nemethi@math.ohio-state.edu;
nemethi@renyi.hu}
\thanks{The first author is supported by the Netherlands Organization for Scientific Research (NWO).
The first three authors are partially
supported by BFM2001-1488-C02-01. The forth
 author is partially supported by NSF grant DMS-0304759.}




\maketitle \pagestyle{myheadings} \markboth{{\normalsize J.
Fern\'andez de Bobadilla, I. Luengo, A. Melle-Hern\'andez, A.
N\'emethi}}{ {\normalsize Rational unicuspidal plane curves}}

\section{Introduction}

It is a very old and interesting open problem to characterize
those collections of embedded topological types of local plane
curve singularities  which may appear as singularities of a
projective plane curve $C$ of degree $d$. (We invite the reader
to consult the articles of Fenske, Flenner, Orevkov, Tono,
Zaidenberg, Yoshihara, and the references therein, for recent
developments.)  The goal of the present article is to give a
complete (topological) classification of those cases when $C$ is
rational and it has a unique singularity which is locally
irreducible (i.e. $C$ is unicuspidal) with one Puiseux pair.

In fact, as a second goal, we also wish to present some of the
techniques which are/might be helpful in such a classification,
and we invite the reader to join us in our effort to produce a
classification  for all the cuspidal rational plane curves. In
fact, this effort also motivates that decision, that in some
cases (in order to have a better understanding of the present
situation), we produce more different arguments for some of the
steps.

In the next paragraph we formulate the main result. We will write
$d$ for the degree of $C$  and $(a,b)$ for the  Puiseux pair of
its cusp,  where $1<a<b$. We denote by $\{\varphi_j\}_{j\geq 0}$
the Fibonacci numbers $\vp_0=0$, $\vp_1=1$,
$\vp_{j+2}=\vp_{j+1}+\vp_j$.

\subsection{Theorem.}\label{maintheorem} {\em The Puiseux pair
$(a,b)$ can be realized by a unicuspidal rational plane curve of degree
$d$ if and only if $(d,a,b)$ appears in the following list.

 \vspace{1mm}

(a) $(a,b)=(d-1,d)$;

(b) $(a,b)=(d/2, 2d-1)$, where $d$ is even;

(c) $(a,b)=(\vp_{j-2}^2,\vp_j^2)$ and $d=\vp_{j-1}^2+1=\vp_{j-2} \vp_{j}$, where $j$ is odd
and $\geq 5$;

(d) $(a,b)=(\vp_{j-2},\vp_{j+2})$ and $d=\vp_{j}$, where $j$ is odd and
$\geq 5$;

(e) $(a,b)=(\vp_4,\vp_8+1)=(3,22)$ and $d=\vp_6=8$;

(f) $(a,b)=(2\vp_4,2\vp_8+1)=(6,43)$ and $d=2\vp_6=16$.}

\vspace{2mm}

All these cases are realizable: (a) e.g. by $\{zy^{d-1}=x^d\}$,
(b) by $\{(zy-x^2)^{d/2}=xy^{d-1}\}$, or by the parametrization
$[z(t):x(t):y(t)]=[1+t^{d-1}:t^{d/2}:t^d]$. The existence of (c)
and (d) is guaranteed by Kashiwara classification \cite{Kash},
Corollary 11.4. These two cases can be realized by a rational
pencil of type $(0,1)$: the generic member of the pencil is (c),
while the special member of the pencil is of type (d) (cf. also
with the last paragraphs of the present article).
 Orevkov in
\cite{Orev} provides a different construction for curves which
realize the case (d) (denoted by him by $C_j$). Similarly, the
cases (e) and (f) are realized by  the sporadic cases $C_4$ and
$C_4^*$ of Orevkov \cite{Orev}.

\subsection{Remarks.}\label{remarkintr} (1) Since $C$ is rational
and its  singular locus $p$ has Milnor number $\mu=(a-1)(b-1)$,
the genus formula says that
\begin{equation*}
(a-1)(b-1)=(d-1)(d-2).\tag{1}\end{equation*} On the other hand,
not any triple $(d,a,b)$ with $(a-1)(b-1)=(d-1)(d-2)$ can be
geometrically realized. E.g., $(5,3,7)$ or $(17,6,49)$ cannot.

(2) There are two integers which coordinate the above
classification. The first one is defined as follows. Let $\pi:X\to
\bp^2$ be the minimal good embedded resolution of $C\subset
\bp^2$, and let $\bar{C}$ be the strict transform of $C$. Clearly,
$(\pi^*C,\cc)=C^2=d^2$, and $\pi^*C =\cc +abE_{-1}+\ldots$ (where
$E_{-1} $ is the unique $-1 $ exceptional curve of $\pi$),
 hence $d^2=\cc^2+ab$.    Using (1), we get:
\begin{equation*} \left\{
\begin{array}{l} a+b=3d-1-\cc^2\\ ab=d^2-\cc^2. \end{array}\right.
\tag{2}
\end{equation*}
Then $\bar{C}^2$ in the above cases is as follows: it is positive
for (a) and (b), it is zero for (c), equals $-1$ for (d), and
$=-2$ for (e) and (f).

(3) The second guiding integer is the logarithmic  Kodaira
dimensions $\bar{\kappa}:=\bar{\kappa}(\bp^2\setminus C)$ (cf.
\cite{F}). Its values are the following (cf. \cite{Orev}):
$-\infty$ for (a)-(d), and 2 for the last two sporadic cases. (In
particular, $\bar{\kappa}$ depends only on the integers
$(d,a,b)$, and it is independent on the analytic type of $C$
which realizes these integers.)

In particular, the above classification shows that
$\bar{\kappa}=-\infty$ if and only if $\bar{C}^2> -2$.

In fact, after we finished the manuscript, we learned from the
introduction of  \cite{Tono} that in \cite{Y} (written in
Japanese) it is proved that for any unicuspidal rational curve
$C$, $\bar{\kappa}=-\infty$ if and only if $\bar{C}^2> -2$. Using
\cite{Y} (i.e. this equivalence), a possible `quick'
classification for $\bar{C}^2> -2$ would run as follows: Since for
all these cases $\bar{\kappa}=-\infty$, we just have to separate
in Kashivara's classification \cite{Kash} those unicuspidal
curves with  exactly one Puiseux pair. Their numerical invariants
$(d,a,b)$ are  exactly those listed in (a)-(d).

On the other hand, this argument probably does not show what is
really behind the classification of this case. Therefore, we
decided to keep the structure of our manuscript, and provide an
independent classification.

Note also that in \cite{BLMN} we list  the complete topological
classification of the cuspidal rational curves with
$\bar{\kappa}<2$. In fact $\bar{\kappa}=0$ cannot occur because
of a result of Tsunoda's \cite{Tsu1}, see also Orevkov's paper
\cite{Orev}. Moreover, Tono in \cite{Tono} provides all the
possible curves $C$ with $\bar{\kappa}=1$: there is no one with
one Puiseux pair.

Hence, in our case,  the remaining part of the classification
corresponds to  $\bar{C}^2\leq -2$, or equivalently, to
$\bar{\kappa}=2$. In general, the classification of this
(`general') case  is the most difficult; and in our case it is not
clear at all at the beginning (and, in fact,  it is rather
surprising) that there are only two (sporadic) cases satisfying
these data.

(4) Let $\alpha=(3+\sqrt{5})/2$ be the root of
$\alpha+\frac{1}{\alpha}=3$. Notice that in family (d) $d/a$ and
$b/d$ asymptotically equals $\alpha$. In fact, for $j$ odd,
$\{\vp_j/\vp_{j-2}\}_j $ are the increasing convergents of the
continued fraction of $\alpha$. Using this, another remarkable
property of the family (d) can be described as follows (cf.
\cite{Orev}, page 658). The convex hull of all the pairs $(m,d)\in
\bz^2$ satisfying $m+1\leq d<\alpha m$ (cf. with the sharp
Orevkov inequality \cite{Orev}, or \ref{Orev})  coincides  with
the convex hull of all pairs $(m,d)$ realizable by rational
unicuspidal curves $C$ (where $d=\deg(C)$ and $m=mult(C,p)$) with
$\bar{\kappa}(\bp^2\setminus C)=-\infty$; moreover, this convex
hull is generated by curves with numerical data (a) and (d).

(5) It is clear that the families (a)-(d) are organized in nice
series of curves. It is less clear from the statement of the
theorem, but rather clear from the proof, that also (e)-(f) form a
`series': they are the only curves with $3d=8a$ (cf. also with the
next remark).

(6) A  hidden massage of the classification (and some of the
steps of the proof) is that there is an intimate relationship
between the semigroup of $\bn$ generated by the elements $a$ and
$b$, and the intervals of type $(\, (l-1)d\, , ld\,]$. The
endpoints $d$ and $3d$ play crucial roles in some of the
arguments. (E.g.,  $\bar{C}^2\leq -2$ if and only if $a+b>3d$;
see also \ref{semigroup}.) This fact is deeply exploited in
\cite{BLMN}. In fact, that paper strongly motivated  the present
manuscript.

(7) The (part of the) proof in section 4 clearly shows the
deficiencies of the known restrictions, bounds which connect the
local data $(a,b)$ with the degree $d$ -- although we list and try
to use a large number of them. On the other hand, the above
classification fits perfectly with the conjectured restriction
proposed by the authors in \cite{BLMN} (valid in a more general
situation), which, in fact,  alone would provide the
classification.

\medskip

The authors thanks Maria Aparecida Ruas 
and Jean-Paul Brasselet for the wonderfull atmosphera 
in the ``VIII\'eme Rencontre Internationale
 de S\~ao Carlos sur les singularit\'es 
re\'elles et complexes au CIRM".
We also thank J. Stevens for pointing us out the full strength 
of the semicontinuity of the spectrum to eliminate some cases in the
classification. 
Finally we thank David Lehavi who wrote for the fourth author  
some computer programs to help in the classification.

\section{Restrictions and bounds}

In the present section we list some general results which impose
some restrictions for the integers $(d,a,b)$. We start with a
trivial one: (1) and (2) clearly imply:

\subsection{The `trivial' bound.}\label{triv} {\em In any
situation $b\geq d$. Moreover,  if $b=d$ then $(a,b)=(d-1,d)$.}

\vspace{2mm}

 If $b>d$ then $a<d-1$. The next result proves a
`gap' for $a$: if $a<d-1$ then $a\leq d/2$ too.

\subsection{Lemma. The `dual curve bound'.}{\label{2} {\em If $b>d$
then $d\geq 2a$ (hence $b>2a$ too).}

\begin{proof} Let $(C,p)$ be the germ of the singular point $p$ of
$C$, and let $\{m_i\}_i$ be the multiplicity sequence of $(C,p)$.
We will use the symbol $^{\vee}$ for the corresponding invariants
of the dual curve $C^{\vee}$ of $C$. By a result of C.T.C. Wall
\cite{Wall} Proposition 7.4.5, the blow ups of the singularities $(C,p)$ and
$(C^{\vee},p^{\vee})$ (where $p^{\vee}$ corresponds to the tangent
cone of $(C,p)$) are equisingular. Assume that $b<2a$. Then
$m_2=b-a$, hence $m_2^\vee=b-a\leq m_1^\vee$. But, according
\cite{Wall}, the intersection multiplicity of the tangent cone of
$(C,p)$ with $C$ at $p$ is $i=m_1+m_1^\vee$, hence $d\geq
i=m_1+m_1^\vee\geq a+b-a=b$, a contradiction. In particular,
$b\geq 2a$. In this case $m_2=a$, hence $m_2^\vee=a$ as well. The
above argument gives: $d\geq i\geq m_1+m_1^\vee\geq 2a$.
\end{proof}

\subsection{The semicontinuity of the spectrum.}\label{4} The very existence
of the curve $C$ shows that the local plane curve singularity
$(C,p)$ is in the deformation of the local plane curve singularity
$(U,0):=(x^d+y^d,0)$ (see e.g. \cite{Dimca} (3.24)). In
particular, we can use the semicontinuity of the spectrum for
this pair \cite{Var1,Var2}. More precisely, this assures that in
any interval $(c,c+1)$, the number of spectral numbers of $(C,p)$
is not larger than the number of spectral numbers of $(U,0)$.
E.g., for the intervals $(-1+l/d,l/d)$
 ($l=2,3,\ldots, d$) one has  the following inequality:
\begin{equation*}
\#\{\frac{i}{a}+\frac{j}{b}<\frac{l}{d}; \ i\geq 1, \ j\geq 1\}
\leq 1+2+\cdots+l-2=\frac{(l-2)(l-1)}{2}.\tag{$SS_l$}\end{equation*}
Notice that the inequality $(SS_d)$ is automatically satisfied
(with equality), since for both singularities the number of spectral
numbers strict smaller than 1 is $(d-1)(d-2)/2$.

\subsection{Example. The inequality $(SS_{d-1})$.}{\label{5} We
denote by $\#_{d-1}$
the number of lattice points at the left hand side of $(SS_{d-1})$.
Since $i/a<(d-1)/d$ and $a<d$, one gets that $1\leq i\leq a-1$. Therefore,
$$\#_{d-1}=
\sum_{i=1}^{a-1} \#\{ j \ : \ 1\leq j < b(\frac{d-1}{d}-\frac{i}{a})\}=
\sum_{i=1}^{a-1} \Big\lceil b-\frac{b}{d}-\frac{ib}{a} \Big\rceil -1,
$$hence
$$\#_{d-1}=(a-1)(b-1)-
\sum_{i=1}^{a-1} \Big\lfloor \frac{b}{d}+\frac{ib}{a} \Big\rfloor.$$
This expression can be computed explicitly. Indeed, since $(a,b)$ is a
lattice point and gcd$(a,b)=1$, one has:
$$\sum_{i=1}^a \Big\lfloor \frac{ib}{a}\Big\rfloor =\frac{(a+1)(b+1)}{2}-a,$$
hence
$$\sum_{i=1}^a \Big\lfloor \frac{ib}{a}+\frac{b}{d}\Big\rfloor =
\frac{(a+1)(b+1)}{2}-a+a\Big\lfloor \frac{b}{d}\Big\rfloor+
\sum_{i=1}^a\Big\lfloor \Big\{\frac{ib}{a}\Big\}+\Big\{\frac{b}{d}\Big\}
\Big\rfloor.$$
Notice that the set $\{ib/a\}$ for $i=1,\ldots,a$ is the same as the set
$r/a$ for $r=0,\ldots,a-1$. Moreover, $r/a+\{b/d\}\geq 1$
if and only if $a-1\geq r\geq \lceil a(1-\{b/d\})\rceil$,
hence the number of possible $r$'s is $\lfloor a\{b/d\}\rfloor$.
Therefore,
$$\sum_{i=1}^a \Big\lfloor \frac{ib}{a}+\frac{b}{d}\Big\rfloor =
\frac{(a+1)(b+1)}{2}-a+\Big\lfloor \frac{ab}{d}\Big\rfloor.$$
Hence
$$\sum_{i=1}^{a-1} \Big\lfloor \frac{ib}{a}+\frac{b}{d}\Big\rfloor =
\frac{(a+1)(b+1)}{2}-a-b-\Big\lfloor \frac{b}{d}\Big\rfloor
+\Big\lfloor \frac{ab}{d}\Big\rfloor,$$
or
$$\#_{d-1}=
\frac{(a-1)(b-1)}{2}+  \Big\lfloor \frac{b}{d}\Big\rfloor
-\Big\lfloor \frac{ab}{d}\Big\rfloor,$$
Then, using (1) and (2), $(SS_{d-1})$ becomes:
\begin{equation*}
\Big\lfloor \frac{b}{d}\Big\rfloor +\Big\lceil
\frac{\cc^2}{d}\Big\rceil \leq 2.\tag{3}\end{equation*}

\subsection{Other examples of $(SS_l)$.}\label{8} $(SS_2)$ is equivalent with
$1/a+1/b\geq 2/d$. This is true automatically, since $1/a+1/b\geq
1/d+1/(d-1)>2/d$. The next  inequality $(SS_3)$ is equivalent with
$2/b+1/a\geq 3/d$, which also is satisfied automatically.

If $b>d$ then  $a+2b>3a+b$ (cf. \ref{2}), hence $(SS_4)$ is
equivalent with the pair of inequalities: $a+2b\geq 4ab/d$ and $
4a+b\geq 4ab/d$. Or, via (2):
$$\min\{3a,b\}\geq  d+1+\frac{d-4}{d}\cc^2.$$
This with an (absolute) lower bound for $\cc^2$ already is interesting:
$3a> d+$const, which has the flavour of the Matsuoka-Sakai inequality $3a>d$
(see \ref{MS}) proved by different methods.

By a similar method as above, one can verify that
 $(SS_{d-2})$ is equivalent with:
\begin{equation*}
\Big\lfloor \frac{2b}{d}\Big\rfloor +\Big\lceil \frac{2\cc^2}{d}\Big\rceil
\leq 5.\end{equation*}
and  $(SS_{d-3})$ is equivalent with:
\begin{equation*}
\Big\lfloor \frac{3b}{d}\Big\rfloor
+\Big\lfloor \frac{3b}{d}-\frac{b}{a}\Big\rfloor
 +\Big\lceil \frac{3\cc^2}{d}\Big\rceil
\leq 8.\end{equation*} In general,  one  expects that the set of
all inequalities $(SS_l)$ is really strong.

\vspace{2mm}

\noindent  The next set of restrictions are provided by {\em
Bogomolov-Miyaoka-Yau type} inequalitities:

\subsection{Matsuoka-Sakai inequality.}\label{MS}  The inequality \cite{MS}
in our case reads as $d<3a$ (valid for any $\bar{\kappa}$).

\subsection{Remark. Orevkov's inequality.}\label{Orev} Orevkov  in \cite{Orev}
obtained  different improved versions of \ref{MS}. Below
$\alpha=(3+\sqrt{5})/2\approx 2.618$ and $\beta=1/\sqrt{5}$.

 \vspace{2mm}

(a) \cite{Orev} Theorem B(a): If $\bar{\kappa}=-\infty$, then
$d<\alpha a$.

(b) \cite{Orev} Theorem B(b): If $\bar{\kappa}=2$, then
$d<\alpha(a+1)-\beta$.

(c) \cite{Orev} (2.2)(4): If $\bar{\kappa}=2$, then
\begin{equation*}-\cc^2\leq -2 +\frac{a}{b}+\frac{b}{a}.\tag{4}
\end{equation*}

\vspace{2mm}

Finally, we end with the following:

\subsection{The `semigroup density property'.}\label{semigroup}
\cite{BLMN} Let $\Gamma$ be the semigroup of $(C,p)$, i.e. the
semigroup (with 0) of $\bn$ generated by the integers $a$ and $b$.
Then for any $0\leq l<d$ the following inequality holds:
$$ \#\Gamma\cap [0,ld]\geq (l+1)(l+2)/2.$$

\begin{proof} It is instructive to sketch the proof for
$l=3$ case:  we wish to prove $\#\Gamma\cap [0,3d]\geq 10$. Recall
that a cubic is determined by 9 parameters. Therefore,
$\#\Gamma\cap [0,3d]\leq 9$ would imply the existence of a cubic
with intersection multiplicity with $C$ at $p$ strict greater
than $3d$, which contradicts B\'ezout's theorem.
\end{proof}

In the classical theory, many `candidates' $(d,a,b)$ were
eliminated by different geometric constrictions using ingenious
Cremona transformations. We will exemplify this in \ref{cre}.

\section{The classification in the case $\cc^2>0$}

\subsection{Theorem.}\label{7} {\em If $\cc^2>0$
then either $(a,b)=(d-1,d)$ or $(a,b)=(d/2,2d-1)$. }

\begin{proof} Since $b\geq d$ (cf. \ref{triv}), by (3) we get that
$\cc^2\leq d$. Clearly, equality holds  if and only if
$(a,b)=(d-1,d)$. Next, assume that $0<\cc^2<d$. Then again   by
$(3)$ one has $\lfloor b/d\rfloor\leq 1$, or $b<2d$. But notice
that $b<2d-1$ would imply (by (1)) that $a>d/2$ which contradicts
\ref{2}. Hence $b=2d-1$.
\end{proof}

\section{Classification in the case $\cc^2\leq -2$}

\subsection{}\label{000}
Our first goal is to prove that $3d\geq 8a$.
For this we apply \ref{semigroup} for $l=3$. Since $a+b>3d$ (cf.
(2)) and $9a>3d$ (cf. \ref{MS}), the needed 10 elements of
$\Gamma\cap [0,3d]$ must be $b,0, a, \ldots, 8a$, hence $8a\leq
3d$. 

\subsection{Corollary.}\label{12} {\em $\bar{\kappa}(\bp^2\setminus C)=2$.}

\begin{proof}
$\bar{\kappa}$ cannot be $-\infty$ because of \ref{Orev}(a);
cannot be 0 because of \cite{Orev}, Theorem B(c) (see also \cite{Tsu1}). Unicuspidal
rational curves with $\bar{\kappa}=1$ are classified by K. Tono
\cite{Tono}, the corresponding splice diagrams are listed in
\cite{BLMN}: there is no example with one Puiseux pair.
\end{proof}

Now, the
classification for $\bar{C}^2\leq -2$ can be finished in two different ways.

\subsection{First proof. Using the computer.}\label{compu}
The first version is based on the inequality \ref{Orev}(b).
Notice that in the case of a geometric realization one must have
$$3\alpha(a+1)-3\beta > 3d\geq 8a,$$
which is true only if $a\leq 44$, (or, by using again
$d<\alpha(a+1)-\beta$), only if $d\leq 117$.
 Hence, we have only to analyze the {\em
finite} family determined by, say,  $d\leq 117$. Then, one can
search with the computer for 3-uples $(d,a,b)$ verifying all the
restrictions considered above.  E.g., we used the conditions
$d\leq 117$, $\gcd(a,b)=1,\, a<d<b,\, d<3a$, $3d\geq 8a$, $2\leq
-\cc^2\leq -2 +\frac{a}{b}+\frac{b}{a}$,
 $b<\alpha (d-1) (d-2)/(d-2\alpha+\beta))+1$,
$(d-\alpha+\beta)/\alpha<a$, and
$(SS_{d-1}),(SS_{d-2}),(SS_{d-3}),(SS_{d-4}), (SS_{4})$.
Using the inequality $3d\geq 8a$ and a similar computation as in
the case of $(SS_{d-1})$, we obtain that $(SS_{d-4})$ is
equivalent with
\begin{equation*}
\Big\lfloor \frac{4b}{d}\Big\rfloor + \Big\lfloor
\frac{4b}{d}-\frac{b}{a}\Big\rfloor + \Big\lceil
\frac{4\cc^2}{d}\Big\rceil \leq 13.\tag{5}\end{equation*} Then
the only triplets satisfying all these are listed below (in the
list appears $(d,a,b;\cc^2)$:

$$
\array{l}
C_1:=(\ 8, 3, 22; -2), \\
C_2:=(11, 4, 31; -3), \\
C_3:=(16, 6, 43; -2), \\
C_4:=(17, 6, 49; -5), \\
C_5:=(19, 7, 52; -3), \\
C_6:=(20, 7, 58; -6).
\endarray
$$
Next, notice that  the curves
 $C_1$ and $C_3$ exist, they are listed in our classification theorem.
The others do not exist:
$C_2$ is eliminated by Orevkov in \cite{Orev}, page 2 (see also \ref{ctwo} (b));
$C_4$ and $C_6$ can be excluded by the semicontinuity property of the
spectrum (applied for all the intervals of type $(l/d, l/d+1)$, $-d<l<d$),
finally $C_5$ can be eliminated by the `nodal
cubic trick', see \ref{ctwo} (a).
(Notice also that $C_2$ and $C_5$ cannot be eliminated by the
semicontinuity property.)


\subsection{Second proof. Resolving diophantic equations.}\label{dio}
Next we show how one can analyse   the case $3d\geq 8a$ (cf.
\ref{000}) by a diophantic equation (for the convenience of the
reader, later we will make more precise the geometry behind this
equation, cf. \ref{155} and \ref{cre}). Our goal is to eliminate
everything excepting $C_1$ and $C_3$, and to emphasize that $C$
exists if and only if $3d=8a$, and $C_1$ and $C_3$ are the only
solutions with $3d=8a$.

Let us write $x:=3d-8a\geq 0$. Then clearly $3|a-x$. Moreover,
$$-\cc^2(a-1)=-(3d-1-a-b)(a-1)=-(x+7a-1-b)(a-1)$$
$$=(b-1)(a-1)-(x+7a-2)(a-1)=(d-1)(d-2)-(x+7a-2)(a-1).$$
Using again $d=(x+8a)/3$ one gets
\begin{equation*}
-9\bar{C}^2(a-1)=x^2+7ax+a^2+9a.
\tag{6}\end{equation*}

\subsection{The case $x=0$.}\label{xzero} (6) implies
the divisibility $a-1|10$. Since one also
has $3|a$,  the only solutions are $a=3$ and $a=6$, corresponding
to $C_1$ and $C_3$ above.

\subsection{Facts.}\label{cbar}  $-\bar{C}^2\leq 7$ {\em and}
$x\leq 5$.

\begin{proof} First we verify $-\bar{C}^2\leq 7$. It is easy to verify
(using (1), (2) and $d/3<a\leq d/2$, cf. \ref{MS} and \ref{2})
that for $6\leq d\leq 10$ this is true. Hence assume that
 $d\geq 11$. Notice that if for some (positive) $k$ one has
$kd\leq -\bar{C}^2<(k+1)d$, then (3) gives $b/d\leq 3+k$. But
$d/a<3$ by \ref{MS}, hence $b/a<3(3+k)$. Using \ref{Orev}(c) one
gets $-\bar{C}^2\leq 3k+7$. Since for $k>0$ and $d\geq 11$ one
has $3k+7 < dk$, one should have $k=0$.

Using this and $x\geq 6$, from (6) we get
$63(a-1)\geq 36+42a+a^2+9a$, which has no solution.
\end{proof}


Now, we consider the above equation (6) for $x\geq 1$.
By \ref{cbar} we only
have to analyse the cases $1\leq x\leq 5$, and eliminate
all the solutions.

\vspace{2mm}

\noindent {\bf The case $x=1$.} In this case one has
$-9(\bar{C}^2+2)(a-1)=(a-1)^2+18$, hence $3|a-1|18$ but
$9\not|a-1$. In particular, $a=4$ or 7
corresponding to $C_2$ and  $C_5$ above.

\noindent {\bf The case $x=2$.} Similarly as above, $a-1|28$ and
$3|a-2$, hence $a-1=4, 7$ or 28. In fact,  if $a=5$ then $d=12$
and $b\not\in \bz$. The next case $(d,a,b;\cc^2) =(22,8,61;-4)$
can be eliminated by $(5)$; the last $(78,29,210;-6)$ by
\ref{Orev}(c).

\noindent {\bf The case $x=3$.} Now $a-1|40$ and $3|a$. The
possible $a$'s are $a=3$ which gives $d=9$ contradicting
\ref{MS}; $a=6$ providing $C_4$; $a=9$ providing $(25,9,70;-5)$
which can be eliminated by $(5)$, and $a=21$ providing
$(57,21,155;-6)$ which is eliminated by \ref{Orev}(c).

\noindent {\bf The case $x=4$.} (6)  has two solutions: $C_6$:
$(20,7,58;-6)$ and $(28,10,79;-6)$, the second one can be
eliminated by (5).

\noindent {\bf The case $x=5$} provides two integral solutions:
 $(23,8,67;-7)$ and $(31,11,88;-7)$.
Both can be eliminated by \ref{Orev}(c).

\vspace{2mm}

We  end this section by the description of the promised geometric
construction (used also
in \cite{Orev} and by E. Artal-Bartolo as well).

\subsection{Lemma. The existence of a specific nodal cubic.}\label{155}
{\em There exists an (unique) irreducible cubic $N\subset \bp^2$
with a node singularity at $p$  such that $N$ and $C$ share the
first seven infinitely near points at $p.$}

\begin{proof}
A cubic is determined by nine parameters. The multiplicity
sequence of $N$ at $p$ should be $[2,1_6].$ Passing through $p$
and having multiplicity 2 provides $3$ conditions. The remaining
six conditions are imposed by the remaining six infinitely near
points. The condition which would imply that
 the singularity $(N,p)$ is a cusp would involve another
equation (the vanishing of the determinant of the quadratic part
at $p$), and the corresponding system of equations would not have
any solution. Similar arguments eliminates other type of
singularities (two smooth branches with contact two, or $(N,p)$
with multiplicity 3). Hence $(N,p)$ is a node.

Next we prove that $N$ cannot be a product of three linear forms.
Indeed, the tangent line $L_0$ of $C$ at $p$ goes just through the
first two infinitely near points because $d<3a$ and $d=L_0 \cdot
C.$  Any other line has less tangency  than $L_0$. This also
shows that $N$ cannot be $L_0\cdot Q$ for some $Q$ (transversal
to $L_0$ at $p$).

The remaining posibility is $N=LQ$ where $Q$ is a smooth conic
and $L$ and $Q$ meets transversally at $p$. Since $Q$ is
determined by five conditions (five infinitely near points) then
$Q$ and $C$ must be tangent and share the seven infinitely near
points at $p.$ In particular by Bezout $2d=Q\cdot C\geq 6a$ which
is in contradiction with $d<3a$, cf. \ref{MS}.
\end{proof}

\subsection{The Cremona transformation associated with the nodal
cubic $N$.}\label{cre}

Consider the nodal cubic $N$ given by  \ref{155}. First we verify that
$C$ and $N$ share exactly the first seven infinitely near points.
Indeed, assume that this is not the case.
If $b\leq 8a$ then the multiplicity sequence of $(C,P)$ is
$[a_7,b-7a,...]$, hence $3d\geq 2a+6a+b-7a=a+b=3d-1-\bar{C}^2>3d$,
a contradiction. If $b>8a$ then the multiplicity sequence of $(C,P)$
is $[a_8,...]$, hence $3d\geq 9a$ which contradicts \ref{MS}.

In particular, the intersection multiplicity of $C$ and $N$ at $P$ is
$8a$. Assume that $C\cap N=\{P,P_1,\ldots,P_r\}$. Notice that at $P_i$
($1\leq i\leq r$) both curves $C$ and $N$ are smooth, let $k_i$ be their
intersection multiplicity at $P_i$. By Bezout's theorem one has
$3d=8a+\sum_ik_i$. We prefer to write $x:=\sum_ik_i$, hence
$3d=8a+x$ (and the notation is compatible with above).

Blow up the common seven infinitely near points. We get seven
irreducible exceptional divisors $\{E_i\}_{i=1}^7$. Let
$\tilde{C}$ and $\tilde{N}$ be the strict transforms of $C$ and
$N$. One has the following intersections: $E_1^2=\cdots=
E_6^2=-2$, \  $ E_7^2=-1$, \ $\tilde{N}^2=-1$, \ $E_1\cdot E_2=E_2\cdot E_3= \cdots
=E_6\cdot E_7=1$,\ $E_1\cdot \tilde{N}=E_7\cdot \tilde{N}=1$.
Also, $\tilde{C} $ intersects $E_7$ (but not the other irreducible
exceptional divisors) at a point $P'$, and the singularity
$(\tilde{C},P')$ has exactly one Puiseux pairs of type
$(b-7a,a)$. The intersection of $\tilde{N}$ with $E_7$ is not
$P'$.

Consider now the curve $\tilde{N}\cup\cup_{i=1}^6E_i$. Clearly, this can be
blown down, and after this modification $\pi$ we get another copy of
$\bp^2$. Let the image of $\tilde{C}$ via this projection be $C'$.
By standard (intersection) argument one gets that the degree
$d'$ of $C'$ is
$$d'=8d-21a \ \ \mbox{(which also satisfies $3d'=8x+a$)}.$$
The curve $C'$ has at most two singular points. One candidate is
the (isomorphic) image of the germ at $P'$ with one Puiseux pair
$(b-7a,a)$. The other is the common image of the points $\{P_i\}$
($1\leq i\leq r$). Clearly, if $x=0$ then this point does not
exist, if $x=1$ then this is a smooth point, but otherwise it is
singular. One can find its embedded resolution graph by blowing
up (for each $i$) $k_i$ times the point $P_i$. Hence, by
A'Campo's formula one can determine its Milnor number, which is
$\mu=7x^2-7x-r+1$ (provided that $x\geq 1$). Since it has $r$
local irreducible  components, the  delta-invariant is
$(7x^2-7x)/2$. Then one can verify that (6) corresponds to the
genus formula of $C'$.

\subsection{Example.}\label{ctwo}
 
(a) Let us start with $(d,a,b)=(19,7,52)$. Then $x=1$, hence $C'$
is again rational and unicuspidal with $(d',a',b')=(5,3,7)$. But
such a curve does not exist because of \ref{7} (one can also check
the classification of rational curves of degree five e.g. in
\cite{Nam}).

(b) Let us consider now the curve $C_2$ above with data
$(d,a,b)=(11,4,31)$. Then $x=1$, hence $C'$ is rational unicuspidal, 
say at $Q_1$, with
$(d',a',b')=(4,3,4)$. Notice that a curve with this triplet may
exists -- although $C_2$ does not. 
The image $\bar{N}$ under the modification $\pi$ of the exceptional curve $E_7$ is a (rational) nodal cubic with a node, say at $Q_2 (\ne Q_1)$. Moreover, $\bar{N}\cdot C'=4 Q_1+ 8 Q_2$. At $Q_1$, $\bar{N}$ is non-singular and with the same tangent as $C'$, and at $Q_2$ the quartic 
$C'$ has intersection multiplicity $7$ with one of the branches of the node of $\bar N$ and 1 with the other.
To show that $C_2$ does not exist we will prove that such configuration of the  rational curves $C'$ and $\bar N$ in $\bp^2$ does not exist.

Choosing affine coordinates we may assume that $C'$ is given by 
the zero locus of $a y^3+a_1 y^3x+a_2 y^2x^2+a_3 y x^3+x^4+a_0 y^4$;
with $a\ne 0.$ In such a case $Q_1=(0,0)$ and its tangent line $L_1=\{y=0\}$ verifies $L_1\cdot C'=4 Q_1$. The curve $C'$ has a parametrization given by $[z(\lambda,t):x(\lambda,t):y(\lambda,t)]=[\lambda^4+a_3t\lambda^3+a_2t^2\lambda^2+a_1t^3\lambda+a_0t^4:-at^3\lambda:-at^4]$.
 
To have $I_{Q_1}(\bar{N}, C')=4$ then $\bar{N}$ must be the zero locus of a polyomial  $y+f_2(x,y)+f_3(x,y)$ (see the parametrization of $C'$), where $f_2(x,y)=m_{1,1}xy+m_{2,0}x^2+m_{0,2}y^2$ and
$f_3(x,y)=n_{1,2} x y^2+n_{2,1}x^2y+n_{3,0} x^3+n_{0,3}y^3$.

Next one imposes that, in the affine plane $\bp^2\setminus L_1=\{y\ne 0\}$, the curves $C'$ and $\bar N$ must meet at only one point $Q_2$
(with intersection multiplicity 8).
The parametrization of $C'$ in this affine chart is $(z,x)=(s^4+a_3s^3+a_2s^2+a_1s+a_0,-as)$ and the equation of $\bar{N}$ is given by
$z^2+f_2(x,1)z+f_3(x,1)=0.$ Imposing to have a solution of the form $(As+B)^8$ gives
$B=a_3A/4$ which means $s=-a_3/4$. We have two possibilities: firstly, if $a_3=0$ then $s=0$
and $Q_2=(z,x)=(a_0,0)$. The solutions are given by
$m_{1,1}=2a_1/a;\,m_{2,0}=-2a_2/a^2;\,m_{0,2}=-(2a_0-a_2^2);\,
n_{1,2}=(-2a_0+a_2^2)a_1/a;\,n_{2,1}=(a_1^2+2a_0a_2-a_2^3)/a^2;\,n_{3,0}=-2a_1a_2/a^3;\,n_{0,3}=(a_0-a_2^2)a_0.$ To have $\bar{N}$ a node at $Q_2$
implies $a_2$ vanishes and therefore $\bar N$ must be a conic which is a contradiction.
 
In the other case, i.e. $a_3\ne 0$ then $s=-a_3/4$ and $Q_2=(z,x)=(z_0,a a_3/4)$.The solutions are given by:
$$
\array{l}
m_{1,1}=(16a_1-a_3^3)/(8a);\qquad \qquad m_{2,0}=(3/4a_3^2-2a_2)/a^2;\\
\\
m_{0,2}=-2a_0+a_2^2+19a_3^4/128-3a_2a_3^2/4;\\
\\
n_{1,2}= -(4096a_0a_1-2048a_1a_2^2-304a_1a_3^4+1536a_1a_2a_3^2-256a_0a_3^3+a_3^7)/(2048a);\\
\\
n_{2,1}=(2^{11}a_0a_2+(2^{5}a_1)^2-2^{10}a_2^3-152a_2a_3^4+768a_2^2a_3^2-128a_1a_3^3-768a_0a_3^2+7a_3^6)/(2^5a)^2;\\
\\
n_{3,0}=(-64a_1a_2+32a_3a_2^2+3a_3^5-20 a_2a_3^3+24 a_1a_3^2)/(32 a^3);\\
\\
n_{0,3}=a_0^2-a_0a_2^2-(19/128)a_0a_3^4+(3/4)a_0a_2a_3^2+(1/65536)a_3^8.
\endarray
$$

In order $\bar{N}$ to have multiplicity two at $Q_2$ one needs
 $a_2=3a_3^2/8$ but this condition also impose 
that the tangent cone of $\bar N$ at $Q_2$ is a double line and therefore $Q_2$
cannot be a node. Hence
this configuration also does not exist.

\section{The case $\cc^2=0,-1$}

In this section we find all the integer solution $(d,a,b)$ of (2)
 with $\cc^2=0,-1$ and we show that all of them can be realized by some
unicuspidal rational plane curve of degree $d$ and Puiseux pair
$(a,b)$. Let $\varphi_j$ be the $i$-th Fibonacci number, that is
$\varphi_0=0, \ \varphi_1=1$ and
$\varphi_{j+2}:=\varphi_{j+1}+\varphi_{j}.$ They share  many
interesting properties, see e.g. \cite{Vaj}. We will use here the
following :
\begin{equation*}
3\varphi_j=\varphi_{j-2}+\varphi_{j+2},\quad \text{and} \quad \varphi_j^2
=(-1)^{j+1}+\varphi_{j-1}\varphi_{j+1}.
\tag{7}\end{equation*}
Let $\Phi=\frac{1+\sqrt{5}}{2}$ be the positive solution of the equation
$\Phi^2-\Phi-1=0.$ For every integer $j>0$ one has :
\begin{equation*}
\Phi^j=\frac{\varphi_{j+1}+\varphi_{j-1}+\varphi_{j}\sqrt{5}}{2}.
\tag{8}\end{equation*}

\subsection{The Pell equation.}\label{17} The system of equations (2)
for $\cc^2=0,-1$ can be transformed (see below) into the {\em
Pell} equation:
\begin{equation*}
x^2-5y^2=-4,\, x,y \in \bz. \tag{9}\end{equation*} Consider the
number field $K=\bq[\sqrt{5}]$ and its ring of integers
$R=\bz[\sqrt{5}]$, which is a UFD. If $\gamma=x+y\sqrt{5}$ is a
solution of (9) then its norm is $N_K(\gamma)=\gamma
\bar{\gamma}=-4.$ Consider $\eta=1+\sqrt{5}$, then $N_K(\eta)=-4$
and $-4$  has a prime decomposition $-4=\eta \bar{\eta}.$ Since
the fundamental unit of $K$ turns out to be $u=2+\sqrt{5}$ and
$\gamma$ is associated either to $\eta$ or $\bar{\eta}$ then
$\gamma$ is either $\pm u^r \eta$ or $\pm \bar{u}^r \bar{\eta}$
(since $\bar{u}=-1/u$) for $r\in\bz.$ Moreover $N_K(u)=-1$ which
implies that $r$ must be even, that is $r=2j$ for $j\in\bz.$ Then
$\eta=2\Phi$ and from the identity $\Phi^2=\Phi+1$ one gets
$\Phi^3=u.$

Thus solutions of (9) are either $\gamma=\pm u^{2j} \eta=\pm 2\Phi^{6j+1}$
or $\gamma=\pm \bar{u}^{2j} \bar{\eta}=\pm 2 \bar{\Phi}^{6j+1}$ with $j\in \bz.$
Using $\Phi \bar{\Phi}=-1$, $\gamma$ is either $\pm 2\Phi^{6j+1},$
$\pm 2\Phi^{6j-1},$ or their conjugates $\pm 2 \bar{\Phi}^{6j+1},\pm 2 \bar{\Phi}^{6j-1}$
with $j\geq 0.$
Using (7) and (8)
the set of solutions of (9) is given by
\begin{enumerate}
\item[(A)] $\pm \left( \varphi_{6j+2}+\varphi_{6j}+\varphi_{6j+1}\sqrt{5} \right),\,\, \text{with} \, j\geq 0, $
\item[(B)] $\pm \left( \varphi_{6j}+\varphi_{6j-2}+\varphi_{6j-1}\sqrt{5} \right),\,\,
\text{with} \, j\geq 0, $
\item[(C)] $\pm \left( \varphi_{6j+2}+\varphi_{6j}-\varphi_{6j+1}\sqrt{5} \right),\,\, \text{with} \, j\geq 0, $
\item[(D)] $\pm \left( \varphi_{6j}+\varphi_{6j-2}-\varphi_{6j-1}\sqrt{5}  \right),\,\,\text{with} \, j\geq 0.$
\end{enumerate}

\subsection{The case $\cc^2=0$.} Since $\gcd(a,b)=1$ and $ab=d^2$ then $a=m^2, b=n^2$  and $d=mn$
for some positive integers $m,n$ with $\gcd(m,n)=1.$ Thus
$a+b=3d-1$ transforms into
\begin{equation*}
m^2+n^2=3mn-1.
\tag{10}\end{equation*}

\subsection{The case $\cc^2=-1$.} The system (2)  provides the equation
\begin{equation*}
a^2+d^2=3ad-1. \tag{11}\end{equation*} Thus, any solution
$(\omega,v)$ of  $\omega^2+v^2=3\omega v-1$ is a solution of
$(2\omega-3v)^2-5v^2=-4$. Hence, with the transformation
$x=2\omega-3v,\  y=v$, one gets the solutions of (9).

{\bf Case A.} If $\gamma=\pm \left(
\varphi_{6j+2}+\varphi_{6j}+\varphi_{6j+1}\sqrt{5} \right),\,\,
j\geq 0,$ is a solution of (9) then $v=\pm \varphi_{6j+1}$ and
$\omega=\pm (\varphi_{6j+2}+\varphi_{6j}+3\varphi_{6j+1})/2 =\pm
\varphi_{6j+3}$ is a solution of (10) and (11) (for the last
equality use (7)). Since $1<a<d$, if $\cc^2=-1,$ then
$a=\varphi_{6j+1},\,d=\varphi_{6j+3}$ and
$b=3d-a=3\varphi_{6j+3}-\varphi_{6j+1}= \varphi_{6j+5}$ for some
$j>0$, by property (7) of Fibonacci numbers. Similarly, if
$\cc^2=0,$ then $\omega$ and $v$ are both either positive or
negative which implies $a=\varphi_{6j+1}^2,\,b=\varphi_{6j+3}^2$
and $d=\omega v=\varphi_{6j+1}\varphi_{6j+3} =\varphi_{6j+2}^2+1$.

{\bf Case B.} If $\gamma=\pm \left(
\varphi_{6j}+\varphi_{6j-2}-\varphi_{6j-1}\sqrt{5} \right),\,\,
j\geq 0, $ is a solution of (9) then $v=\pm (-\varphi_{6j-1})$ and
$\omega=\pm (\varphi_{6j}+\varphi_{6j-2}-3\varphi_{6j-1})/2=\pm
(-\varphi_{6j-3})$ is a solution of (10) and (11). In the case
$\cc^2=-1,$ one gets  $a=\varphi_{6j-3},\,d=\varphi_{6j-1}$ and
$b=3d-a=3\varphi_{6j-1}-\varphi_{6j-3} =\varphi_{6j+1}$ with
$j>0.$ If $\cc^2=0,$ then $\omega$ and $v$ are both either
positive or negative which implies
$a=\varphi_{6j-3}^2,\,b=\varphi_{6j-1}^2$ and $d=\omega
v=\varphi_{6j-1}\varphi_{6j-3} =\varphi_{6j-2}^2+1$ with $j>0.$

{\bf Case C.} If $\gamma=\pm \left( \varphi_{6j+2}+\varphi_{6j}-\varphi_{6j+1}\sqrt{5} \right),\,\, j\geq 0, $
is a solution of (9)
then $v=\pm (-\varphi_{6j+1})$ and
$\omega=\pm (\varphi_{6j+2}+\varphi_{6j}-3\varphi_{6j+1})/2=\pm (-\varphi_{6j-1})$
is a solution of (10) and (11).
If $\cc^2=-1,$ then $a=\varphi_{6j-1},\,d=\varphi_{6j+1}$ and $b=\varphi_{6j+3}$
with $j>0.$
If $\cc^2=0,$ then $\omega$ and $v$ are both either positive or negative
which implies $a=\varphi_{6j-1}^2,\,b=\varphi_{6j+1}^2$ and $d=\varphi_{6j-1}\varphi_{6j+1} =\varphi_{6j}^2+1$
with $j>0.$

{\bf Case D.} Any solution in this case is included in the
previous cases.

\vspace{1mm}

 Hence,  we determined  all the possible integer  solutions.

\subsection{Theorem. Classification for $\cc^2=-1$}\label{19} {\em If $\cc^2=-1$
then $(a,b)=(\varphi_{j-2},\varphi_{j+2})$ and $d=\varphi_j$, with $j$ odd $\geq 5.$ For every such $j$ there exists
a unicuspidal rational plane curve of degree
with such invariants.}

\subsection{Theorem. Classification for $\cc^2=0$}\label{20} {\em If $\cc^2=0$
then $(a,b)=(\varphi_{j-2}^2,\varphi_{j}^2)$ and $d=\varphi_{j-1}^2+1,$ with $j$ odd $\geq 5.$
For every such $j$ there exists
a unicuspidal rational plane curve with such invariants.}

\medskip

We only need to provide the equations of the curves. We will rely
on \cite{Kash},  Corollary 11.4. Let $(x,y)$ be a system of affine
coordinates in $\bp^2$ and consider
\begin{equation*}
P_{-1}=y-x^2,\,\, Q_{-1}=y,\,\,P_0=(y-x^2)^2-2xy^2(y-x^2)+y^5,
\end{equation*}
\begin{equation*}
Q_0=y-x^2,\,\, G=xy-x^3-y^3,\,\, Q_s=P_{s-1}, \,\,  P_s=\left(G^{\varphi_{2s+1}}+Q_s^3\right)/Q_{s-1}.
\end{equation*}

Then $P_s$ is a polynomial in $x$ and $y$ of degree
$\varphi_{2s+3}$ and defines a  rational unicuspidal curve whose
unique singularity $p$ has exactly one characteristic pair of type
$(a,b)=(\varphi_{2s+1},\varphi_{2s+5}).$ The curves $P_s=0$ and
$Q_s=0$ only meet at $p$. The rational pencil with only one base
point determined by the rational function
$R_s=(P_s)^{\varphi_{2s+1}}/(Q_s)^{\varphi_{2s+3}}$ has only two
special fibres $P_s=0$ and $Q_s=0$,  and the other fibres
are rational unicuspidal plane curves of degree
$\varphi_{2s+3}\varphi_{2s+1}=\varphi_{2s+2}^2+1.$ The singularity
of a generic fiber has  one characteristic pair
$(a,b)=(\varphi_{2s+1}^2,\varphi_{2s+3}^2).$

\end{document}